\documentclass[graybox]{svmult}

\usepackage{type1cm}        % activate if the above 3 fonts are
\usepackage{makeidx}         % allows index generation
\usepackage{graphicx}        % standard LaTeX graphics tool
\usepackage{multicol}        % used for the two-column index
\usepackage[bottom]{footmisc}% places footnotes at page bottom

\usepackage{newtxtext}       % 
\usepackage[varvw]{newtxmath}       % selects Times Roman as basic font

\usepackage{enumerate}

\makeindex             % used for the subject index
\begin{document}

\title*{Error estimates and variance reduction for nonequilibrium stochastic dynamics}
\titlerunning{Nonequilibrium stochastic dynamics} 
\author{Gabriel Stoltz}
\institute{Gabriel Stoltz \at CERMICS, Ecole des Ponts \& Matherials team-project, Inria Paris \email{gabriel.stoltz@enpc.fr}}
\maketitle

\abstract{Equilibrium properties in statistical physics are obtained by computing averages with respect to Boltzmann--Gibbs measures, sampled in practice using ergodic dynamics such as the Langevin dynamics. Some quantities however cannot be computed by simply sampling the Boltzmann--Gibbs measure, in particular transport coefficients, which relate the current of some physical quantity of interest to the forcing needed to induce it. For instance, a temperature difference induces an energy current, the proportionality factor between these two quantities being the thermal conductivity. From an abstract point of view, transport coefficients can also be considered as some form of sensitivity analysis with respect to an added forcing to the baseline dynamics. There are various numerical techniques to estimate transport coefficients, which all suffer from large errors, in particular large statistical errors. This contribution reviews the most popular methods, namely the Green--Kubo approach where the transport coefficient is expressed as some time-integrated correlation function, and the approach based on longtime averages of the stochastic dynamics perturbed by an external driving (so-called nonequilibrium molecular dynamics). In each case, the various sources of errors are made precise, in particular the bias related to the time discretization of the underlying continuous dynamics, and the variance of the associated Monte Carlo estimators. Some recent alternative techniques to estimate transport coefficients are also discussed. }

%-------------- intro to MD ---------------------
\section{An introduction to computational statistical physics}
\label{sec:intro}

Statistical physics offers a way to pass from a microscopic description of systems to a macroscopic one. Its numerical realization is known as ``molecular dynamics'', which is a simulation method used and developed over the past 70~years; see~\cite{BaCiGr20} for a historical perspective, and~\cite{AlTi17,FrSm02,Tu10} for reference textbooks in the physics literature. Molecular dynamics has two major aims. First, it can be used as a {\it numerical microscope}, which allows to perform ``computer'' experiments. This was the initial motivation for simulations at the atomic scale: physical theories were tested on computers. Another major aim of molecular simulation is to compute macroscopic quantities or thermodynamic properties, typically through averages of some functionals of the system. This allows to obtain a \emph{quantitative} information on a system, instead of resorting to approximate theories constructed for simplified models. Molecular dynamics can therefore be seen as a tool to explore the links between the microscopic and macroscopic properties of a material, allowing to address modeling questions such as ``Which microscopic ingredients are necessary (and which are not) to observe a given macroscopic behavior?''

I start by describing in this section how to compute static properties, such as the pressure at a given temperature and density. I then turn to the main topic in Section~\ref{sec:LR}, namely the computation of transport coefficients, such as the thermal conductivity or the shear viscosity. The numerical analyses of methods in molecular dynamics to estimate these transport coefficients are provided in Section~\ref{sec:error}. I highlight in particular why the computation of these transport coefficients is much more difficult than the computation of static properties. Current open questions and research tracks to address this point are given in Section~\ref{sec:ccl}.

\subsection{Passing from a microscopic to a macroscopic description}

Physical systems are described at the microscopic level by their positions~$q \in \mathcal{D} = \mathbb{R}^d$ or~$(L\mathbb{T})^d$ (with~$\mathbb{T} = \mathbb{R}\backslash \mathbb{Z}$ the one-dimensional torus), and momenta $p \in \mathbb{R}^{d}$. The associated phase-space is denoted by $\mathcal{E} = \mathcal{D} \times \mathbb{R}^d$. The interactions between the particles are taken into account through a potential energy function~$V:\mathcal{D}\to \mathbb{R}$, depending on the positions~$q$ only. The total energy of the system is given by the Hamiltonian
\begin{equation}
  \label{intro:Hamiltonian}
  H(q,p) = V(q) + \frac{1}{2} \, p^\top M^{-1} p,
\end{equation}
where~$M \in \mathbb{R}^{d \times d}$ is a positive definite matrix (called the mass matrix, most usually a diagonal matrix). The macroscopic state of a system is described, within the framework of statistical physics, by a probability measure~$\mu$ on the phase space $\mathcal{E} = {\cal D} \times \mathbb{R}^{d}$. Macroscopic properties are obtained via averages of an observable~$\varphi$ with respect to this measure:
\begin{equation}
  \label{intro:averages}
\mathbb{E}_\mu(\varphi) = \int_{\mathcal{E}} \varphi(q,p) \, \mu(dq \, dp).  
\end{equation}

In many physical situations, systems in contact with some energy thermostat are considered, rather than isolated systems with a fixed energy. In this case, the energy of the system fluctuates. A typical way to describe such systems is to consider their microscopic configurations to be distributed according to the so-called \emph{canonical measure}, which is the following Boltzmann--Gibbs probability measure: 
\begin{equation}
  \label{intro:canonical}
  \mu(dq \, dp) = Z_\mu^{-1} \exp (-\beta H(q,p)) \, dq \, dp,
\end{equation}
where $\beta>0$ is proportional to the inverse of the temperature, and~$Z_\mu$ is the normalization constant. The canonical measure is of tensorized form 
\[
\mu(dq \, dp) = \nu(dq) \, \kappa(dp),
\]
where 
\begin{equation}
  \label{intro:density}
  \nu(dq) = Z_\nu^{-1} \textrm{e}^{-\beta V(q) } \, d q,
  \qquad
  Z_\nu = \int_{\cal D} {\rm e}^{-\beta V(q)}\, dq,
\end{equation}
and~$\kappa$ is a Gaussian measure with covariance~$M/\beta$. Therefore, sampling configurations $(q,p)$ according to the canonical measure~$\mu(dq\,dp)$ can be performed by independently sampling positions according to~$\nu(dq)$ and momenta according to~$\kappa(dp)$. Since it is straightforward to sample from~$\kappa$, the actual issue is to sample from~$\nu$. 

\subsection{Computing average properties with Langevin dynamics}

The main mathematical challenge in computing averages such as~\eqref{intro:averages} is the very high dimensionality of the integral under consideration, which prevents the use of standard quadrature methods. In practice, the only realistic option is to rely on ergodic averages over trajectories of dynamics leaving the probability measure~$\mu$ invariant. In this context, the average of some observable~$\varphi \in L^1(\mu)$ is approximated as
\begin{equation}
  \label{eq:ergodic_avg}
  \int_\mathcal{E} \varphi \, d\mu = \lim_{t \to +\infty} \widehat{\varphi}_t, \qquad \widehat{\varphi}_t = \frac1t \int_0^t \varphi(q_s,p_s) \, ds.
\end{equation}
I focus here on Langevin dynamics
\begin{equation}
  \label{eq:Langevin}
  \left \{ \begin{aligned}
    d q_t & = M^{-1} p_t \, dt, \\
    d p_t & = - \nabla V (q_t) \, dt - \gamma \, M^{-1} p_t \, dt + \sqrt{\frac{2\gamma}{\beta}} \, dW_t,
    \end{aligned} \right.
\end{equation}
where~$W_t$ is a standard $d$-dimensional Brownian motion, and~$\gamma > 0$ the magnitude of the friction term. This dynamics can be seen as a stochastic perturbation of the Hamiltonian dynamics. In order to mathematically study the properties of the Langevin dynamics, and understand for instance why the prefactor~$\sqrt{2\gamma/\beta}$ in front of the Brownian motion is the right one, it is useful to introduce some operators, in particular the evolution semigroup
\[
\left(\mathrm{e}^{t \mathcal{L}}\varphi\right)(q,p) = \mathbb{E}\left[\varphi(q_t,p_t) \, \Big| (q_0,p_0)=(q,p)\right],
\]
and its generator 
\begin{equation}
  \label{eq:generator_Langevin}
  \mathcal{L} = \mathcal{L}_\mathrm{ham} + \gamma \mathcal{L}_\mathrm{FD},
\end{equation}
where the Hamiltonian and fluctuation/dissipation parts of the evolution are encoded by the generators
\[
\mathcal{L}_\mathrm{ham} = p^\top M^{-1}\nabla_q - \nabla V^\top \nabla_p,
\qquad
\mathcal{L}_\mathrm{FD} = - p^\top M^{-1}\nabla_p + \frac1\beta \Delta_p.
\]
The invariance of the probability measure~$\mu$ is then characterized by the property (denoting by~$C^\infty_{\rm c}(\mathcal{E})$ the vector space of smooth functions with compact support)
\begin{equation}
  \label{eq:invariance_mu_by_L}
  \forall \varphi \in C^\infty_{\rm c}(\mathcal{E}), \qquad \int_\mathcal{E} \mathcal{L} \varphi \, d\mu = 0,
\end{equation}
which is equivalent to~$\mathcal{L}^\dagger \mu = 0$, the latter equality being in the sense of distributions, with~$\mathcal{L}^\dagger$ the adjoint of~$\mathcal{L}$ on~$L^2(\mathcal{E})$.

From a functional analytical viewpoint, it is in fact more convenient to work in~$L^2(\mu)$, because the operators appearing in the decomposition~\eqref{eq:generator_Langevin} are respectively antisymmetric and symmetric when considered on this functional space. Moreover, $L^2(\mu)$ is the natural functional space to consider to state central limit theorems and quantify the statistical error. To make these considerations precise, we denote by~$A^*$ the adjoint of a closed operator~$A$ on~$L^2(\mu)$. Then, 
\[
\mathcal{L}^* = -\mathcal{L}_\mathrm{ham} + \gamma \mathcal{L}_\mathrm{FD},
\qquad
\mathcal{L}_\mathrm{FD} = -\frac1\beta \sum_{i=1}^d \partial_{p_i}^* \partial_{p_i},
\qquad
\mathcal{L}_\mathrm{ham} = \frac{1}{\beta} \sum_{i=1}^d \partial_{p_i}^* \partial_{q_i} - \partial_{q_i}^* \partial_{p_i}.
\]
Indeed, a simple computation gives  
\[
\int_{\mathcal{D}} \left(\partial_{q_i} \varphi\right) \phi \, d\mu = -\int_{\mathcal{D}} \varphi \left(\partial_{q_i} \phi\right) d\mu - \int_{\mathcal{D}} \varphi \phi \, \partial_{q_i} \left(Z_\nu^{-1} \mathrm{e}^{-\beta V} \right) d\kappa,
\]
so that $\partial_{q_i}^* = -\partial_{q_i} + \beta \partial_{q_i} V$. A similar computation gives $\partial_{p_i}^* = -\partial_{p_i} + \beta (M^{-1} p)_i$. In this functional framework, the invariance of a probability measure~$\rho$ with density~$f$ with respect to~$\mu$ can be reformulated as 
\[
\mathcal{L}^* f = 0.
\]
For Langevin dynamics, a simple argument based on the (anti)symmetry of the two terms in~\eqref{eq:generator_Langevin} allows to prove that~$f=\mathbf{1}$ is the unique solution in~$L^2(\mu)$ of the previous equation; see~\cite[Proposition~15]{Vi09}. Indeed, one first proves that~$\nabla_p f = 0$ by multiplying the equation by~$f$ and integrating with respect to~$\mu$; therefore, $f$ does not depend on~$p$ and the equation reduces to~$(L_\mathrm{ham}f)(q,p) = p^\top M^{-1} \nabla_q f(q) = 0$, hence~$f$ does not depend on~$q$ either, and is therefore constant. 

\subsection{Ergodicity results for Langevin dynamics}
\label{sec:ergodicity_Langevin}

The almost-sure convergence of the ergodic averages~$\widehat{\varphi}_t$ in~\eqref{eq:ergodic_avg} follows from he results of~\cite{Kl87} since the stochastic dynamics preserves a probability measure with positive density, and its generator is hypoelliptic~\cite{Ho67}. The asymptotic variance of these ergodic averages allows to quantify the statistical error:
\begin{equation}
  \label{eq:asymptotic_variance}
  \begin{aligned}
  \lim_{t \to +\infty} t\mathrm{Var}\left[\widehat{\varphi}_t^2\right] & = 2 \int_0^{+\infty} \mathbb{E}_0\left[\Pi_0\varphi(q_t,p_t)\Pi_0\varphi(q_0,p_0)\right] \, dt\\
  & = 2 \int_\mathcal{E} \int_0^{+\infty}\!\! \left(\mathrm{e}^{t \mathcal{L}}\Pi_0 \varphi\right) \Pi_0 \varphi \, dt \, d\mu \\
  & = 2 \int_\mathcal{E} \left(-\mathcal{L}^{-1}\Pi_0 \varphi\right) \Pi_0 \varphi \, d\mu
  \end{aligned}
\end{equation}
where
\begin{equation}
  \label{eq:def_Pi_mu}
  \Pi_0 \varphi = \varphi - \mathbb{E}_\mu(\varphi),
\end{equation}
and where the following operator equality was used:
\begin{equation}
  \label{eq:inv_by_int}
  -\mathcal{L}^{-1} = \int_0^{+\infty} \mathrm{e}^{t \mathcal{L}} \, dt
\end{equation}
on the Hilbert space
\begin{equation}
  \label{eq:L_2_0}
  L^2_0(\mu) = \Pi_0 L^2(\mu) = \left\{ \varphi \in L^2(\mu) \, \left| \int_\mathcal{E} \varphi \, d\mu = 0 \right. \right\}.
\end{equation}
Let us emphasize that~$-\mathcal{L}^{-1}$ can be defined only for functions of average~0 with respect to~$\mu$ since functions $\varphi = \mathcal{L}\Phi$ necessarily have average~0 with respect to~$\mu$ by~\eqref{eq:invariance_mu_by_L}. The definition~\eqref{eq:inv_by_int} is legitimate when the operator norm of the semigroup~$\mathrm{e}^{t\mathcal{L}}$ decays sufficiently fast, for instance exponentially. In fact, in such a setting, the Poisson equation 
\[
-\mathcal{L} \Phi = \Pi_0 \varphi = \varphi - \int_\mathcal{E} \varphi \, d\mu
\]
has a unique solution in $L^2_0(\mu)$, and so a central limit theorem holds~\cite{Bh82}. There are various techniques for obtaining the exponential convergence of the semigroup in Banach subspaces of~$L^2_0(\mu)$. Let us list a few prominent examples, referring to the introduction of~\cite{BeFaLeSt20} for a more extensive review:
\begin{itemize}
\item a first approach is based on Lyapunov techniques~\cite{Wu01,MaStHi02,Re06,HaMa11}, which rely on convergence estimates in the Banach space
  \begin{equation}
    \label{eq:B_infty_K}
    B^\infty_\mathcal{K}(\mathcal{E}) = \left\{ \varphi \, \textrm{measurable}, \, \sup \left|\frac{\varphi}{\mathcal{K}}\right| <+\infty \right\},
  \end{equation}
  where~$\mathcal{K} : \mathcal{E} \to [1,+\infty]$ is a Lyapunov function, \emph{i.e.} a function such that~$\mathcal{L} \mathcal{K} \leq -a \mathcal{K} + b$ for some constants~$a>0$ and~$b \in \mathbb{R}$. These approaches were historically the first ones to be devised, and also work for dynamics for which the force is not the gradient of a potential. One limitation of Lyapunov techniques is that they are usually not very quantitative, because it is difficult to obtain sharp minorization conditions.
\item the hypocoercive framework~$H^1(\mu)$, popularized by the monograph~\cite{Vi09}, builds on various previous works where commutators of the building blocks of the generator~$\mathcal{L}$ are used in order to retrieve some dissipation in the~$q$ variables from the dissipation in the~$p$ variables~\cite{Ta02,EcHa03,HeNi04}. Convergence results in~$H^1(\mu)$ can then be transferred to~$L^2(\mu)$ after hypoelliptic regularization~\cite{He07};
\item  a more direct route to proving the convergence in $L^2$ was first proposed in~\cite{He06}, then extended in~\cite{DoMoSc09,DoMoSc15}. This method can be applied to Fokker--Planck and Boltzmann-type operators. It is based on a modification of the $L^2$ scalar product with some regularization operator. This more direct approach makes it even easier to quantify convergence rates; see~\cite{DoKlMoSc13,GrSt16} for studies on the dependence of parameters such as friction in Langevin dynamics; 
\item it was recently shown in~\cite{AlArMoNo19,CaLuWa19,Br21} how to directly obtain convergence in~$L_0^2(\mu)$, based on a space-time Poincar\'e inequality involving the operator~$\partial_t - \mathcal{L}_{\rm ham}$;
\item finally, exponential convergence can also be obtained from purely probabilistic arguments based on a clever coupling between two realizations of the stochastic differential equation~\eqref{eq:Langevin}, see~\cite{EbGuZi19}. 
\end{itemize}

%-----------------------
\subsection{Variance reduction methods in equilibrium molecular dynamics}
\label{sec:variance_reduction_eq}

The computation of time averages~\eqref{eq:ergodic_avg} can be difficult in practice because of the possibly high variance of the estimator~$\widehat{\varphi}_t$. This is often due to the metastability of the dynamics, \emph{i.e.} the fact the target measure to sample has several modes of high probability which are separated by low probability regions. For nonequilibrium systems, we will see in Section~\ref{sec:NEMD} that the variance of the time averages of interest is large even for systems which are not metastable.

Two standard variance reduction techniques for equilibrium dynamics such as~\eqref{eq:Langevin} rely on the explicit knowledge of the invariant measure of the system: 
\begin{itemize}
\item in importance sampling methods, the potential energy function~$V$ is replaced by~$V+\widetilde{V}$, with~$\widetilde{V}$ chosen so that the resulting dynamics (\emph{e.g.} \eqref{eq:Langevin} with force~$-\nabla (V+\widetilde{V})$ instead of~$-\nabla V$) is less metastable. The (marginal distribution in the position variable of the) target measure to sample then has a density proportional to~$\mathrm{e}^{-\beta (V+\widetilde{V})(q)}$. Averages with respect to the baseline canonical measure can be recovered by unbiasing the dynamics using weights~$\mathrm{e}^{\beta \widetilde{V}(q_t)}$, relying on estimators such as
  \[
  \widehat{\varphi}_{\widetilde{V},t} = \frac{\displaystyle \int_0^t \varphi(q_s,p_s) \, \mathrm{e}^{\beta \widetilde{V}(q_s)} \, ds}{\displaystyle \int_0^t \mathrm{e}^{\beta \widetilde{V}(q_s)} \, ds}.
  \]
  A key point of the method is that the invariant probability measure of the dynamics~\eqref{eq:Langevin} with force~$-\nabla (V+\widetilde{V})$ can be easily related to the invariant probability measure of the dynamics~\eqref{eq:Langevin} with force~$-\nabla V$.
\item stratification is a way of decomposing a difficult sampling problem into several easier ones. Ideally, the phase space should be decomposed into the collection of all metastable states, corresponding to local minima of the potential energy function, and these regions should be independently sampled. The local averages in each region should then be reweighted according to the canonical weight of the region itself. There are two major ways to make this idea practical, depending on whether the considered regions overlap -- but in both situations one needs the explicit expression of the measure restricted to the region of interest.

  When there is some overlap between the regions, bridge sampling methods such as MBAR~\cite{ShCh08} can be used. The method was theoretically studied in several works in statistics~\cite{Ge94,MeWo96,KoMcMeNiTa03,Ta04}.

  Non-overlapping regions can also be constructed as the level sets of some real-valued function of the configuration of the system. In this case, the sampling is performed by constraining the dynamics on the iso-surfaces corresponding to various values of the level-set function, and varying the values of the constraint in order to sample the full phase space. This method is known as thermodynamic integration, with a reconstruction performed by computing the free energy; see for instance~\cite[Chapter~3]{LeRoSt10} and references therein. A key element in the approach is that the invariant probability measure of the dynamics restricted to the submanifold is the conditioning of the invariant probability measure of the unconstrained dynamics. 
\end{itemize}
In contrast, control variate methods do not require the expression of the invariant measure of the dynamics, and can therefore be used for equilibrium or nonequilibrium systems, as discussed in Section~\ref{sec:current_perspectives}.

%---------------------- LR -------------------
\section{Definition of transport coefficients}
\label{sec:LR}

\begin{figure}[b]
\sidecaption
\includegraphics[width=0.95\textwidth]{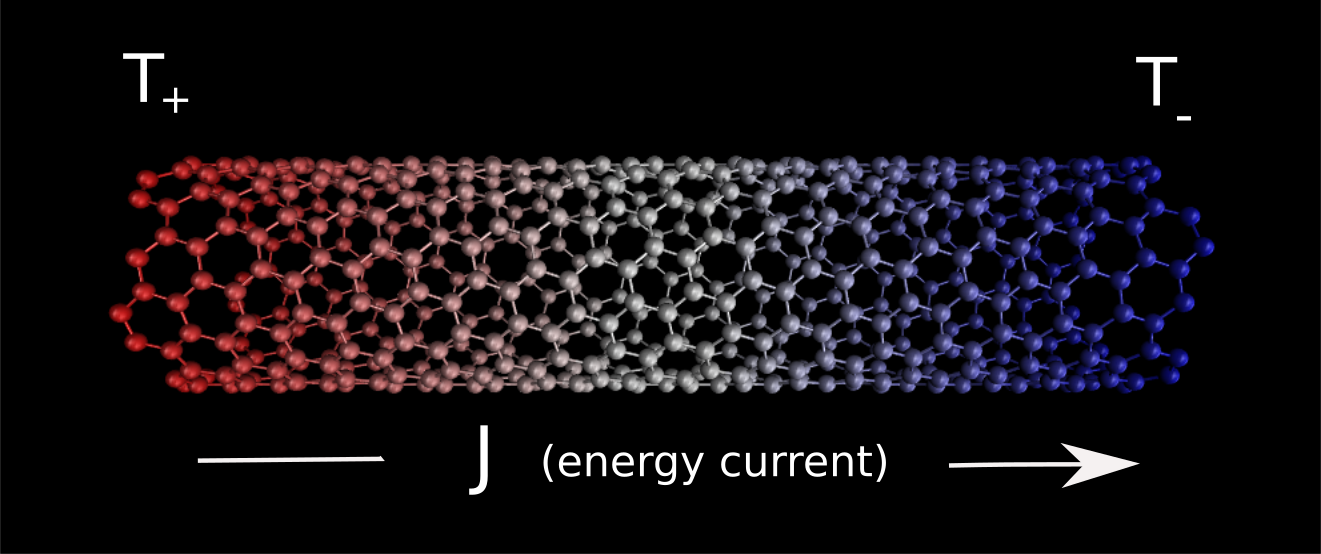}
\caption{Thermal transport in a carbon nanotube.}
\label{fig:CNT}
\end{figure}

\begin{figure}[b]
  \sidecaption
\includegraphics[width=0.8\textwidth]{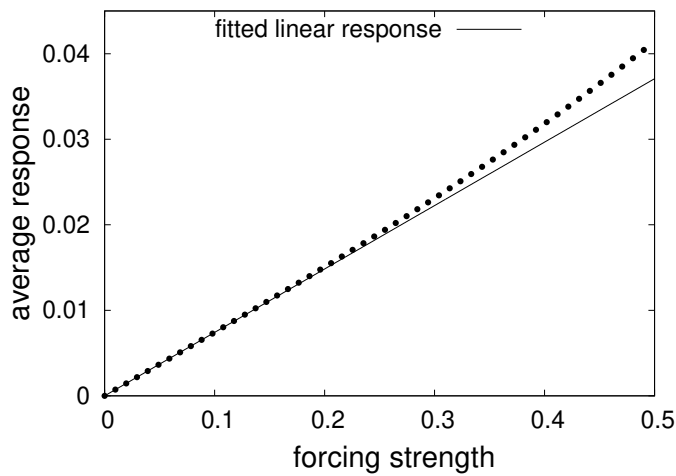}
\caption{Expected steady-state response as a function of the forcing magnitude.}
\label{fig:expected_response}
\end{figure}

At the macroscopic level, transport coefficients relate an external forcing acting on the system (electric field, temperature gradient, velocity field, etc) to an average response expressed through some steady-state flux (of charged particles, energy, momentum, etc). A typical example is Fourier's law
\[
J(x) = - \kappa(x) \nabla T(x),
\]
which relates the energy flux~$J$ to the temperature gradient~$\nabla T$ using the thermal conductivity~$\kappa$ (which is a matrix in general); see Figure~\ref{fig:CNT} for an illustration. At the microscopic level, as discussed in Section~\ref{sec:noneq_dyn}, systems subjected to an external forcing in their steady state are modeled by adding a perturbation of magnitude~$\eta \in \mathbb{R}$ to the reference, equilibrium dynamics. This corresponds to the so-called ``nonequilibrium molecular dynamics'' method~\cite{CiKaSe05,EvMo08,Tu10}. It is observed that, in general, the response~$\mathbb{E}_\eta(R)$ of the system, as encoded by the steady-state average of the physical observable~$R$ of interest, which has average~0 at equilibrium (such as the velocity, the energy flux, etc), is proportional to the magnitude~$\eta$ of the forcing for small values of~$|\eta|$:
\[
\mathbb{E}_\eta(R) \approx \alpha \eta.
\]
This corresponds to the so-called linear response regime, made precise in Section~\ref{sec:LR_pert}. By definition, transport coefficients are the proportionality constants~$\alpha$ relating the response to the forcing; see Figure~\ref{fig:expected_response}. I explain and motivate in Section~\ref{sec:error} why it is difficult to numerically estimate transport coefficients, as Monte Carlo methods suffer from a slow convergence due to a large noise to signal ratio. Long computational times are therefore required to estimate them.

\subsection{Nonequilibrium dynamics and their steady states}
\label{sec:noneq_dyn}

To make precise the mathematical ideas underpinning nonequilibrium molecular dynamics, I consider a specific example, namely the computation of the mobility of fluid particles. This system is modeled as particles enclosed in a periodic cubic box of side~$L>0$, so that positions belong to~$\mathcal{D} = (L\mathbb{T})^d$. A constant force $F \in \mathbb{R}^d$ of norm~1, with an intensity~$\eta \in \mathbb{R}$, is applied to one or several particles in the system. Let us emphasize that a constant force is not the gradient of a smooth periodic function, so that the perturbation to the dynamics is indeed of nongradient type. The corresponding perturbed Langevin dynamics reads
\begin{equation}
  \label{eq:noneq_Langevin}
  \left \{ \begin{aligned}
    dq_t & = M^{-1} p_t \, dt, \\
    dp_t & = \Big( -\nabla V(q_t) + \eta F \Big)dt - \gamma M^{-1} p_t \, dt + \sqrt{\frac{2\gamma}{\beta}} \, dW_t.
  \end{aligned} \right.
\end{equation}
The quantity of interest is the velocity gained in the direction~$F$ by the particles experiencing the forcing, which corresponds to the steady-state average of the response function~$R(q,p) = F^\top M^{-1}p$. The existence and uniqueness of an invariant probability measure for nonequilibrium dynamics such as~\eqref{eq:noneq_Langevin} can be proved using Lyapunov techniques~\cite{Re06,HaMa11}, as made precise in~\cite[Section~5]{LeSt16}. 

In contrast to equilibrium dynamics such as~\eqref{eq:Langevin}, there is no simple analytical expression of the invariant probability measure for the nonequilibrium dynamics~\eqref{eq:noneq_Langevin}, which prevents for instance resorting to Metropolis-type algorithms, and also to the standard variance reduction techniques mentioned in Section~\ref{sec:variance_reduction_eq} (conditioning, biasing, ...). To further discuss properties of the invariant probability measure, it is convenient to introduce the generator~$\mathcal{L} + \eta \widetilde{\mathcal{L}}$ with~$\widetilde{\mathcal{L}} = F^\top \nabla_p$. The invariant probability measure is then a solution to the stationary Fokker--Planck equation
\begin{equation}
  \label{eq:FP_eta_dagger}
  \left( \mathcal{L} + \eta \widetilde{\mathcal{L}} \right)^\dagger \psi_\eta = 0.
\end{equation}
Results of hypoellipticity~\cite{Ho67,Re06} show that the unique invariant probability measure admits a smooth density~$\psi_\eta$ with respect to the Lebesgue measure. One key feature of this invariant measure is that it depends non-locally on the potential~$V$, in contrast to invariant measures~\eqref{intro:canonical} of equilibrium dynamics, for which a local modification of~$V$ leads to a local modification of the density (up to an overall normalization factor). More generally, changes in~$V$ lead to nontrivial modifications in~$\psi_\eta$, which prevents in particular resorting to variance reduction methods such as importance sampling. To illustrate this point, consider the one dimensional overdamped Langevin dynamics~$d q_t =  ( -V'(q_t) + \eta  ) \, d t + \sqrt{2} \, d W_t$ on the torus~$\mathbb{T}$, with~$V \neq 0$ (otherwise the unique invariant probability measure is the uniform measure, whatever the value of~$\eta \in \mathbb{R}$). The density of the unique invariant probability measure satisfies the stationary Fokker--Planck equation
\[
\frac{d}{dq}\left( (V' - \eta) \psi_\eta + \psi_\eta' \right) = 0.
\]
There is therefore a constant~$a \in \mathbb{R}$ such that 
\[
\frac{d}{dq}\left[\psi_\eta(q) \mathrm{e}^{V(q)-\eta q}\right] = a \mathrm{e}^{V(q)-\eta q},
\]
and so, by integration between~$q$ and~$q+1$, 
\[
\psi_\eta(q) \mathrm{e}^{V(q)-\eta q} (\mathrm{e}^{-\eta}-1) = a \int_q^{q+1} \mathrm{e}^{V(z)-\eta z} \, dz.
\]
The constant~$a$ is finally determined by the condition that the periodic function~$\psi_\eta$ has integral~1 over the unit torus: 
\[
\psi_\eta(q) = Z_\eta^{-1} \int_\mathbb{T} \,\mathrm{e}^{V(q+y)-V(q)-\eta y} \, d y.
\]
The latter expression highlights that, when~$\eta \neq 0$, a modification to~$V$ at a given point modifies~$\psi_\eta$ everywhere in a nontrivial way.

\subsection{The linear response regime}
\label{sec:LR_pert}

Although the expression of~$\psi_\eta$ is unknown in general, it can be expanded in powers of~$\eta$ for small values of~$|\eta|$. This corresponds to the so-called perturbative regime. From a mathematical viewpoint, this amounts to writing the invariant probability measure as~$\psi_\eta = f_\eta \mu$ with~$f_\eta = 1 + \mathrm{O}(\eta)$. In order to identify the equation satisfied by~$f_\eta$, one writes
\[
\forall \varphi \in C^\infty_{\rm c}(\mathcal{E}), \qquad 0 = \int_\mathcal{E} \left[ \left(\mathcal{L}+\eta \widetilde{\mathcal{L}}\right)\varphi \right] f_\eta \, d\mu = \int_\mathcal{E} \varphi \left[ \left(\mathcal{L}+\eta \widetilde{\mathcal{L}}\right)^*f_\eta \right] d\mu,
\]
where we recall the notation~$A^*$ introduced after~\eqref{eq:invariance_mu_by_L}. The Fokker--Planck equation~\eqref{eq:FP_eta_dagger} can then be rewritten as
\[
\left( \mathcal{L}+\eta \widetilde{\mathcal{L}} \right)^*f_\eta = 0.
\]
By identifying powers of~$\eta$ in the latter equation (recalling the definition~\eqref{eq:def_Pi_mu} for~$\Pi_0$), one obtains  
\begin{equation}
  \label{eq:f_eta_expansion}
f_\eta = 1 + \eta \mathfrak{f}_1 + \eta^2 \mathfrak{f}_2 + \dots, \qquad -\mathcal{L}^* \mathfrak{f}_1 = \widetilde{\mathcal{L}}^* \mathbf{1} = \Pi_0 \widetilde{\mathcal{L}}^* \mathbf{1} = S.
%\Big(\mathrm{Id} + \eta (\widetilde{\mathcal{L}}\Pi_0\mathcal{L}^{-1}\Pi_0)^* \Big)^{-1} \mathbf{1} = \left(1 + \sum_{n=1}^{+\infty} (-\eta)^n \left[ (\widetilde{\mathcal{L}}\Pi_0\mathcal{L}^{-1}\Pi_0)^*\right]^n \right)\mathbf{1}
\end{equation}
The function
\begin{equation}
  \label{eq:conjugate_response}
  S = \widetilde{\mathcal{L}}^* \mathbf{1}
\end{equation}
is called the conjugate response. For the running example~\eqref{eq:noneq_Langevin}, simple computations reveal that~$\widetilde{\mathcal{L}}^* = -\widetilde{\mathcal{L}} + \beta F^\top M^{-1}p$, so that
\[
S(q,p) = \beta F^\top M^{-1}p.
\]

The linear response of a property of interest~$R \in L^2_0(\mu) = \Pi_0 L^2(\mu)$ can finally be formulated as follows:
\begin{equation}
  \label{eq:LR_formula}
  \begin{aligned}
    \alpha & = \lim_{\eta \to 0} \frac{\mathbb{E}_\eta(R)}{\eta}
    = \int_\mathcal{E} R \mathfrak{f}_1 \, d\mu = \int_\mathcal{E} R \left[\left(-\mathcal{L}^*\right)^{-1}S\right] d\mu = \int_\mathcal{E} \left(-\mathcal{L}^{-1}R\right) S \, d\mu \\
    & =  \int_0^{+\infty} \left[ \int_\mathcal{E} \left(\mathrm{e}^{t \mathcal{L}} R\right) S \, d\mu \right] dt = \int_0^{+\infty} \!\! \mathbb{E}_0\Big(R(q_t,p_t)S(q_0,p_0) \Big)dt,
  \end{aligned}
\end{equation}
where the expectation in the last equality is taken over initial conditions sampled according to~$\mu$, and realizations of the equilibrium dynamics~\eqref{eq:Langevin}. The last equality corresponds to the so-called Green--Kubo formula, further discussed in Section~\ref{sec:GK}. It can be extended to other nonequilibrium forcings when~$S = \widetilde{\mathcal{L}}^* \mathbf{1} \in L^2_0(\mu)$.

To conclude this section, let us summarize the steps required in practice to define transport coefficients from the perspective of computational statistical physics:
\begin{enumerate}[(1)]
\item The first step is to identify a reference dynamics and the physically relevant response function (which depends on the transport coefficient which is targetted).
\item One should next construct a physically meaningful perturbation. These perturbations can be on the drift of the dynamics, or on the fluctuation term (as for thermal transport, where various parts of the system are maintained at different temperatures). Moreover, the forcing can be either a bulk one, applied to all particles in the system, or boundary driven, acting only on the particles at the boundary or close to it.
\item Finally, the transport coefficient $\alpha$ of interest (thermal conductivity, shear viscosity, ...) can be obtained using numerical methods based on the linear response formula~\eqref{eq:LR_formula}, as made precise in Section~\ref{sec:error}. For the running example~\eqref{eq:noneq_Langevin}, the transport coefficient of interest is the mobility matrix~$D$, defined with the response function~$R(q,p) = F^\top M^{-1}p$ as
\[
\lim_{\eta \to 0} \frac{\mathbb{E}_\eta\left(F^\top M^{-1}p \right)}{\eta} = \beta \, F^\top D F,
\]
where
\[ %begin{equation}
  %\label{eq:D}
  D = \int_0^{+\infty} \!\!\mathbb{E}_0\Big( (M^{-1}p_t) \otimes (M^{-1}p_0) \Big) \, dt.
\] %end{equation}
%\item Non physical forcings giving same transport coefficient (``synthetic'')
\end{enumerate}

%-------------------- error estimates -------------------
\section{Error estimates for the computation of transport coefficients}
\label{sec:error}

I summarize in this section the numerical analysis for the two most prominent numerical methods to estimate transport coefficients, namely nonequilibrium molecular dynamics (Section~\ref{sec:NEMD}), which is based on the first equality in~\eqref{eq:LR_formula}, and the Green--Kubo method (Section~\ref{sec:GK}), based on the last equality in~\eqref{eq:LR_formula}. In both cases, I make precise the various sources of bias and variance in the estimators at hand. I finally briefly discuss in Section~\ref{sec:comparison_NEMD_GK} advantages and limitations of both methods, motivating that there is currently no clear cut preference on which one to choose.

%-------------------- NEMD ------------------------
\subsection{Nonequilibrium molecular dynamics}
\label{sec:NEMD}

The principle of nonequilibrium molecular dynamics is to approximate the limit~$\eta \to 0$ in the first equality~\eqref{eq:LR_formula} by considering a small but finite value of~$\eta$ and replacing the expectation with respect to the steady-state probability measure by a time average. More precisely, recalling the nonequilibrium dynamics~\eqref{eq:noneq_Langevin} (and explicitly indicating the dependence on~$\eta$ by superscripts~$\eta$), 
\[
\left \{ \begin{aligned}
  dq_t^\eta & = M^{-1} p_t^\eta \, dt, \\
  dp_t^\eta & = \Big( -\nabla V(q_t^\eta) + \eta F \Big)dt - \gamma M^{-1} p_t^\eta \, dt + \sqrt{\frac{2\gamma}{\beta}} \, dW_t,
\end{aligned} \right.
\]
an estimator of the linear response is, for an observable~$R$ with average~0 under the equilibrium measure:
\begin{equation}
  \label{eq:as_limit_NEMD}
  \widehat{A}_{\eta,t} = \frac{1}{\eta t}\int_0^t R(q_s^\eta,p_s^\eta) \, ds \xrightarrow[t\to+\infty]{\mathrm{a.s.}} \alpha_\eta := \frac1\eta \int_\mathcal{E} R \, f_\eta \, d\mu = \alpha + \mathrm{O}(\eta).
\end{equation}
The various sources of error for the estimator~$\widehat{A}_{\eta,t}$, made precise below, are the following: 
\begin{enumerate}[(i)]
\item A statistical error with asymptotic variance $\mathrm{O}(t^{-1}\eta^{-2})$, much larger than the usual asymptotic variance of order~$1/t$ associated with standard time averages which are not divided by a factor~$\eta$.
\item A bias of order~$\mathrm{O}(\eta)$ due to the fact that~$\eta \neq 0$, as made apparent on the right hand side of~\eqref{eq:as_limit_NEMD}.
\item A bias arising from the finiteness of the integration time~$t$ in the estimator~$\widehat{A}_{\eta,t}$.
\item A bias arising from the discretization in time when implementing nonequilibrium dynamics in computer simulations.
\end{enumerate}
Let us already emphasize that there is a balance between taking~$\eta$ sufficiently small in order to limit the bias~$\alpha_\eta - \alpha = \mathrm{O}(\eta)$, and~$\eta$ not too small so that the asymptotic variance controlling the magnitude of the statistical error is not too large.

\subsubsection{Analysis of the variance and the finite integration time bias}
Let us first announce the two important results concerning the errors related to the finiteness of the integration time. First, the statistical error is dictated by a Central Limit Theorem: it can be shown that
\begin{equation}
  \label{eq:CLT_NEMD}
\sqrt{t} \left(\widehat{A}_{\eta,t} - \alpha_\eta \right) \xrightarrow[t \to +\infty]{\mathrm{law}} \mathcal{N}\left(0,\frac{\sigma_{R,\eta}^2}{\eta^2}\right),
\end{equation}
where
\[
\sigma_{R,\eta}^2 = \sigma_{R,0}^2 + \mathrm{O}(\eta),
\qquad
\sigma_{R,0}^2 = 2  \int_0^{+\infty} \mathbb{E}_0\left[R(q_t,p_t) R(q_0,p_0)\right] \, dt.
\]
This quantifies the fact that the statistical error~$\widehat{A}_{\eta,t} - \alpha_\eta$ is of order~$1/(\eta \sqrt{t})$.

Note that the asymptotic variance is, at dominant order in~$\eta$, the one for the equilibrium dynamics. This result shows that long simulation times of order~$t \sim \eta^{-2}$ are required in order to have an asymptotic variance of order~1. Balancing the bias of order~$\eta$ arising from~$\alpha_\eta - \alpha = \mathrm{O}(\eta)$ and the statistical error requires~$t \sim \eta^{-3}$. Concerning the finite time integration bias, the following estimate holds:
\begin{equation}
  \label{eq:FT_bias}
  \left| \mathbb{E}\left(\widehat{A}_{\eta,t}\right) - \alpha_\eta \right| \leq \frac{K}{\eta t}.
\end{equation}
The bias due to $t < +\infty$ is therefore of order~$\mathrm{O}\left(t^{-1}\eta^{-1}\right)$, typically smaller than the statistical error.

The key tool for proving these results is the Poisson equation
\[
-\left(\mathcal{L}+\eta\widetilde{\mathcal{L}}\right) \mathscr{R}_\eta = R - \int_\mathcal{E} R f_\eta \, d\mu.
\]
A simple computation based on It\^o calculus gives
\begin{equation}
  \label{eq:A_eta_t_Ito}
  \begin{aligned}
    \widehat{A}_{\eta,t} - \frac1\eta \!\int_{\mathcal{E}} \!R f_\eta \, d\mu & = \frac{\mathscr{R}_\eta(q_0^\eta,p_0^\eta) - \mathscr{R}_\eta(q_t^\eta,p_t^\eta)}{\eta t} \\
    & \quad + \frac{\sqrt{2\gamma}}{\eta t\sqrt{\beta}} \int_0^t \!\!\nabla_p \mathscr{R}_\eta(q_s^\eta,p_s^\eta)^\top dW_s.
  \end{aligned}
\end{equation}
The limit~\eqref{eq:CLT_NEMD} then follows from a Central Limit Theorem for martingales, while~\eqref{eq:FT_bias} is obtained by taking expectations in the above equality. Of course, some care has to be taken in order to make these computations rigorous, as one needs to ensure that the solution~$\mathscr{R}_\eta$ to the Poisson is sufficiently regular, and has good integrability properties; see~\cite{SpSt22} for details.

\subsubsection{Analysis of the timestep discretization bias}

In order to discuss the bias arising from discretizing in time the dynamics, one first needs to make precise the numerical schemes used in practice to integrate dynamics such as~\eqref{eq:noneq_Langevin}. These numerical schemes correspond to a Markov chain characterized by the evolution operator
\[
(P_{\Delta t} \varphi)(q,p) = \mathbb{E}\Big( \varphi\left(q^{n+1},p^{n+1}\right)\Big| (q^n,p^n) = (q,p)\Big).
\]
For Langevin like dynamics, it is convenient to rely on splitting methods, as extensively studied in~\cite{LeMa13,LeMa15,LeMaSt16}. The generator of the dynamics~\eqref{eq:noneq_Langevin} is decomposed into the three following parts, which can all be analytically integrated:
\[
A = M^{-1} p \cdot \nabla_q, 
\quad
B_\eta = \left(-\nabla V(q) + \eta F \right) \cdot \nabla_p,
\quad
C = -M^{-1} p \cdot \nabla_p + \beta^{-1} \Delta_p.
\]
First- and second-order splittings are determined by the order in which the various operators~$A,B,C$ are integrated. For example, the scheme with evolution operator~$P_{\Delta t}^{B_\eta,A,\gamma C}$ corresponds to
\begin{equation}
\label{eq:Langevin_splitting}
\left\{ \begin{aligned}
\widetilde{p}^{n+1} & = p^n + {\Delta t} \left( -\nabla V(q^{n}) + \eta F\right), \\
q^{n+1} & = q^n + {\Delta t} \, M^{-1} \widetilde{p}^{n+1}, \\
p^{n+1} & = \rho_{\Delta t} \widetilde{p}^{n+1} + \sqrt{\beta^{-1}(1-\rho_{\Delta t}^2)M} \, G^n, 
\end{aligned} \right.
\end{equation}
where $G^n$ are i.i.d. standard Gaussian random variables and $\rho_{\Delta t} = \exp(-\gamma M^{-1} {\Delta t})$. It can be shown that the so-constructed numerical schemes admit a unique invariant probability measure~$\mu_{\gamma,\eta,{\Delta t}}$, see~\cite{MaStHi02,LeMaSt16,DuEnMoSt21}.

The results established in~\cite{LeMaSt16} provide error estimates à la Talay--Tubaro~\cite{TaTu90}, and show that there exists an integer~$a$, larger than or equal to the weak order of the method, such that 
\[
\int_\mathcal{E} R \, d{\mu}_{\gamma,\eta,{\Delta t}} = \int_\mathcal{E} R \Big(1+ \eta f_{0,1,\gamma} + {\Delta t}^a f_{a,0,\gamma} + \eta {\Delta t}^a f_{a,1,\gamma} \Big) d{\mu} + r_{\varphi,\gamma,\eta,{\Delta t}},
\]
with~$f_{0,1,\gamma} = \mathfrak{f}_1$ is the perturbation~\eqref{eq:f_eta_expansion} of the invariant measure arising from the nonequilibrium forcing, and where the remainder is compatible with linear response:
\[
\left|r_{\varphi,\gamma,\eta,{\Delta t}}\right| \leq K(\eta^2 + {\Delta t}^{a+1}), 
\qquad 
\left|r_{\varphi,\gamma,\eta,{\Delta t}} - r_{\varphi,\gamma,0,{\Delta t}}\right| \leq K \eta (\eta + {\Delta t}^{a+1}).
\]
A corollary of this equality is the following error estimate on the numerically computed mobility:
\[
\begin{aligned}
\alpha_{{\Delta t}} & = \lim_{\eta \to 0} \frac{1}{\eta} \left(\int_\mathcal{E} F^\top M^{-1} p \, \mu_{\gamma,\eta,{\Delta t}}(d{q}\,d{p}) - \int_\mathcal{E}  F^\top M^{-1} p \, \mu_{\gamma,0,{\Delta t}}(d{q}\,d{p}) \right) \\
& = \alpha + {\Delta t}^a \int_\mathcal{E}  F^\top M^{-1} p  \, f_{a,1,\gamma}(q,p) \, \mu(dq \, dp) + {\Delta t}^{a+1} r_{\gamma,{\Delta t}}.
\end{aligned}
\]
It is further possible to state results which hold uniformly in the overdamped limit~$\gamma \to +\infty$; see~\cite{LeMaSt16}.

\begin{figure}[b]
\sidecaption
\includegraphics[width=0.49\textwidth]{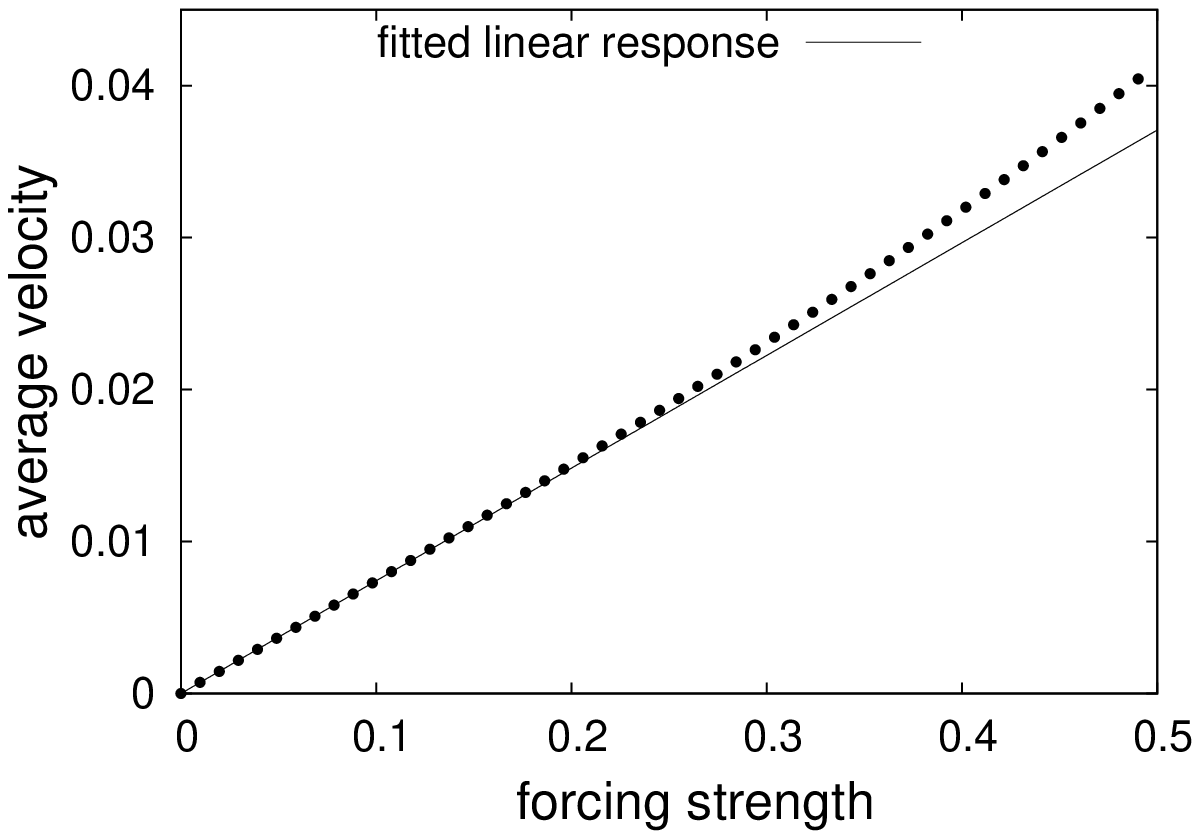}
\includegraphics[width=0.49\textwidth]{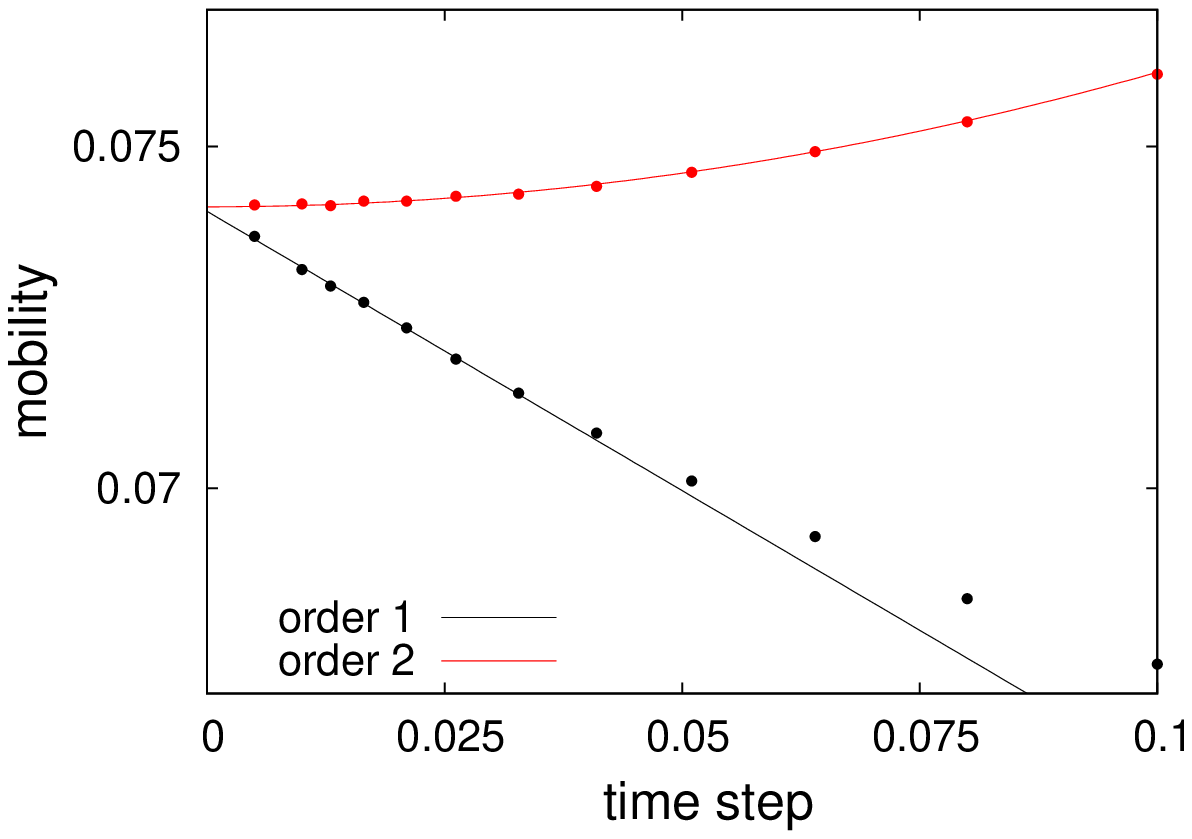}
\caption{Left: Linear response of the average velocity as a function of $\eta$ for the scheme associated with $P_{\Delta t}^{\gamma C, B_\eta,A,B_\eta, \gamma C}$ and ${\Delta t} = 0.01, \gamma = 1$. Right: Scaling of the mobility for the first order scheme $P_{\Delta t}^{A,B_\eta,\gamma C}$ and the second order scheme $P_{\Delta t}^{\gamma C, B_\eta,A,B_\eta, \gamma C}$.}
\label{fig:noneq}
\end{figure}

The error estimates discussed in this section are illustrated for a two-dimensional periodic potential in Figure~\ref{fig:noneq} (see~\cite{LeMaSt16} for details on the numerical computations). The left plot illustrates the linearity of the response for small forcings. The right plot summarizes the estimates for the transport coefficient, obtained from the slopes of the linear responses for various choices of the timestep. Symmetric second order splittings allow to estimate the mobility with an error of order~$\Delta t^2$, whereas first order splitting schemes exhibit much larger errors of order~$\Delta t$.

%------------- GK --------------------
\subsection{Green--Kubo formulas}
\label{sec:GK}

I discuss in this section how to approximate transport coefficients based on the Green--Kubo formula
\begin{equation}
  \label{eq:alpha_GK}
  \alpha = \int_0^{+\infty} \mathbb{E}_0\Big(R(q_t,p_t)S(q_0,p_0) \Big)dt.
\end{equation}
Although the latter expression only requires to refer to equilibrium dynamics, it still implicitly depends on the nonequilibrium forcing through the conjugated response~$S$ defined in~\eqref{eq:conjugate_response}. 

A natural estimator for the right hand side of~\eqref{eq:alpha_GK} is obtained by first truncating the upper bound of the integration time to some time~$T < +\infty$, then replacing the expectation over realizations of the equilibrium dynamics~\eqref{eq:Langevin} by an average over a finite number~$K$ of such realizations, which leads to
\begin{equation}
  \label{eq:natural_GK}
  \widehat{A}_{K,T} = \frac1K \sum_{k=1}^K \int_0^T R(q_t^k,p_t^k)S(q_0^k,p_0^k)\, dt.
\end{equation}
On top of that, the dynamics needs to be discretized using a timestep~$\Delta t>0$, and the time integrals need to be approximated by quadrature formulas. In summary, three sources of errors can be distinguished here: 
\begin{enumerate}[(i)]
\item The truncation of the integration time, usually small given the exponential convergence of the semigroup~$\mathrm{e}^{t \mathcal{L}}$ discussed in Section~\ref{sec:ergodicity_Langevin}. 
\item The statistical error in the estimation of the correlation terms, which increases substantially with the time lag (as discussed in~\cite{DeGr17}). 
\item The bias arising from the finiteness of the timestep and the quadrature formulas for the time integral.
\end{enumerate}
Before embarking on a more precise discussion of these errors, let us mention that it may be beneficial to consider estimators going beyond the natural estimator~\eqref{eq:natural_GK}, using in particular Fourier approaches and time series analysis~\cite{ErMaBa17}, and/or importance sampling on trajectory space~\cite{DoHaKe17}.

\subsubsection{Truncation of time and statistical error}

The truncation bias is small due to the generic exponential decay of correlations: when~$\|\mathrm{e}^{t \mathcal{L}} \varphi\|_{L^2(\mu)} \leq C\|\varphi\|_{L^2(\mu)} \mathrm{e}^{-\kappa t}$ for~$\varphi \in L^2_0(\mu)$ (see Section~\ref{sec:ergodicity_Langevin} for conditions ensuring this), the truncation bias can be upper bounded by
\[
\left|\mathbb{E}\left(\widehat{A}_{K,T}\right)-\alpha\right| = \left| \int_T^{+\infty} \left(\int_\mathcal{E} \left(\mathrm{e}^{t \mathcal{L}}R\right)S\, d\mu \right)dt \right| \leq \frac{C}{\kappa}\|R\|_{L^2(\mu)} \|S\|_{L^2(\mu)} \mathrm{e}^{-\kappa T}.
\]
This suggests to take~$T$ large enough in order to have a small upper bound. On the other hand, the statistical error of~$\widehat{A}_{K,T}$ typically increases with~$T$, so that a trade-off has to be found in the choice of~$T$. More precisely, the analysis in~\cite{PlStWa22} shows that there exists~$C \in \mathbb{R}_+$ such that
\[
\forall T \geq 1, \qquad \mathrm{Var}\left(\widehat{A}_{K,T}\right) \leq C \frac{T}{K}.
\]
The proof of this bound is based on the following equality, obtained as for~\eqref{eq:A_eta_t_Ito} by It\^o calculus on the solution~$\mathscr{R}$ to the Poisson equation~$-\mathcal{L}\mathscr{R} = R \in L^2_0(\mu)$:
\begin{equation}
  \label{eq:GK_Ito}
  \int_0^T R(q_t,p_t) \, dt = \mathscr{R}(q_0,p_0) - \mathscr{R}(q_t,p_t) + \sqrt{\frac{2\gamma}{\beta}} \int_0^T \nabla_p \mathscr{R}(q_t,p_t)^\top dW_t.
\end{equation}
This suggests indeed that~$\widehat{A}_{K,T}$, which is an average over~$K$ realizations of the random variable
\[
\int_0^T S(q_0,p_0) R(q_t,p_t) \, dt \approx S(q_0,p_0) \sqrt{\frac{2\gamma}{\beta}} \int_0^T \nabla_p \mathscr{R}(q_t,p_t)^\top dW_t,
\]
should have a variance scaling as~$T$ since the martingale term on the right-hand side of the previous equality has variance of order~$T$.

\subsubsection{Timestep bias for Green--Kubo formulas}

I discuss in this section the bias on the Green--Kubo formula arising from discretizing the underlying stochastic process. The result is presented for general diffusion processes satisfying certain technical conditions, namely:
\begin{itemize}
\item The bias on the invariant probability measure~$\mu_{\Delta t}$ of the numerical scheme is of order~${\Delta t}^a$ (which is guaranteed under some ergodicity conditions as soon as the numerical scheme of weak order~$a$, see for instance~\cite{LeSt16}, as well as~\cite{AbViZy15,LeMaSt16} for Langevin dynamics).
\item The evolution operator can be expanded as~$P_{\Delta t} = \mathrm{Id} + {\Delta t} \mathcal{L} + {\Delta t}^2 L_2 + \dots + {\Delta t}^{a} L_a + \dots$, with a control of the remainder terms.
\item Uniform-in-$\Delta t$ convergence holds true. It is meant by this that~$P_{{\Delta t}}^{\lceil T/{\Delta t} \rceil}$ decays exponentially with~$T$ in an appropriate functional setting, typically as an operator on the weighted spaces~\eqref{eq:B_infty_K}. Such a decay property is a consequence of the existence of a Lyapunov function, and a minorization conditions uniform in the timestep for initial conditions~$X_0$ in compact sets: for any compact set~$\mathscr{C}$, there exist~$\rho > 0$ and a probability measure~$m$ such that
  \[
  \forall X_0 \in \mathscr{C}, \qquad P_{{\Delta t}}^{\lceil T/{\Delta t} \rceil}\big(X_0,dX\big) \geq \rho \, m(dX).
  \]
  The latter estimate has been established in~\cite{DuEnMoSt21} for splitting schemes of Langevin dynamics. 
\end{itemize}
Under these conditions, the following Riemann-like formula was proved for Langevin dynamics in~\cite{LeMaSt16}, and generalized in~\cite{LeSt16} to more general diffusion processes: For two functions~$R,S \in L^2(\mu)$ with average~0 with respect to~$\mu$,  
\begin{equation}
  \label{eq:GK_dt}
\int_0^{+\infty} \mathbb{E}_0 \Big( R(X_t) S(X_0) \Big) d{t} = {\Delta t} \sum_{n=0}^{+\infty} \mathbb{E}_{\Delta t} \left(\widetilde{R}_{{\Delta t}}\left(X^{n}\right)S(X^0)\right) + \mathrm{O}({\Delta t}^a),
\end{equation}
with
\[
\widetilde{R}_{{\Delta t}} = \widetilde{\mathcal{R}}_{\Delta t} - \mu_{\Delta t}(\widetilde{\mathcal{R}}_{\Delta t}),
\qquad
\widetilde{\mathcal{R}}_{\Delta t} = \Big(\mathrm{Id} + {\Delta t} \,L_2 \mathcal{L}^{-1} + \dots + {\Delta t}^{a-1} L_a\mathcal{L}^{-1} \Big)R.
\]
Simple computations reveal that the right hand side of~\eqref{eq:GK_dt} reduces to a discretization of the time integral with a trapezoidal rule for schemes of weak order~2. Let us also emphasize that an important side result of the above estimate is that the asymptotic variance for numerical schemes coincides at dominant order in~$\Delta t$ with the asymptotic variance of the underlying continuous process, upon considering the dynamics on the same timescale (\emph{i.e.} comparing the variance of trajectory averages computed over times~$\tau$ for the continuous dynamics, and for a number $\tau/\Delta t$ timesteps for the discrete process).

\begin{figure}
\sidecaption
\includegraphics[width=0.8\textwidth]{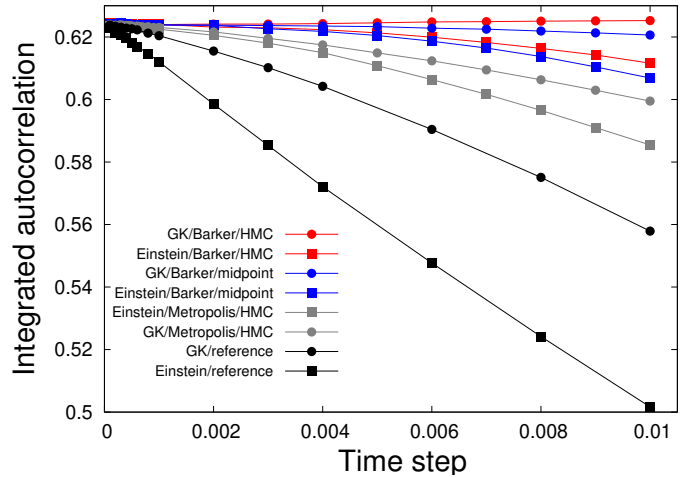}
\caption{Numerical estimates of the integrated autocorrelation of~$V'$ as a function of the timestep~$\Delta t$. See~\cite{FaSt15} for further precision.}
\label{fig:GK}
\end{figure}

The error estimates highlighted in this section have been numerically verified in~\cite{FaSt15}, on a one dimensional system in a periodic potential, governed by Langevin dynamics in the overdamped limit; see Figure~\ref{fig:GK}. The reference numerical method, which corresponds to a numerical scheme of weak order~1 obtained by correcting a Euler--Maruyama scheme with a Metropolis procedure, allows to approximate time integrated correlation functions with an error of order~$\Delta t$. This bias can be reduced by considering more elaborate numerical schemes, and other Metropolis rules, the analysis being based in all cases on~\eqref{eq:GK_dt}.

%---------
\subsection{Advantages and limitations of nonequilibrium molecular dynamics and Green--Kubo formulas}
\label{sec:comparison_NEMD_GK}

The main interest of Green--Kubo formulas is that it suffices to run the equilibrium dynamics, and that integrated correlation functions in~\eqref{eq:alpha_GK} can be computed for various choices of~$R,S$ for the corresponding realizations of the equilibrium dynamics. This allows to potentially compute several transport coefficients with a single simulation. In contrast, a different nonequilibrium perturbation has to be considered for every transport coefficient in the nonequilibrium molecular dynamics approach. Moreover, one needs in principle to check the linearity of the response in nonequilibrium simulations, typically by computing the response for several values of the forcing magnitude, which further increases the computational cost (although these computations can be parallelized in a straightforward way).

However, the postprocessing of Green--Kubo simulations is less straightforward than for nonequilibrium simulations, for which plain time averages~\eqref{eq:as_limit_NEMD} are considered. Moreover, it is possible to increase the magnitude of the forcings in nonequilibrium simulations in order to enhance the response to measure; whereas the correlation functions appearing in the Green--Kubo cannot be enhanced in such a way and may therefore require more computational time to emerge out of the statistical noise.

In the end, the choice of the method to use mostly depends on the habits and experience of the practitioners and the capabilities of the software used. It is a good practice to compare the results of Green--Kubo and nonequilibrium approaches to make sure that the estimated transport coefficients agree. It is also fair to say that various extensions or modifications can or should be used in order to improve the computational workflow of each approach; see Section~\ref{sec:current_perspectives}.

%---------------- perspectives ---------------
\section{Extensions and perspectives}
\label{sec:ccl}

I discuss in this final section alternatives to the two standard approaches to compute transport coefficients described in Section~\ref{sec:error}. The method presented in Section~\ref{sec:CLR_MP} is based on a fluctuation identity different from the Green--Kubo identity. There are however many possible other ways to estimate transport coefficients, which have not been explored in a systematic manner. I give some possible research directions to this end in Section~\ref{sec:current_perspectives}.

\subsection{An example of alternative fluctuation formulas}
\label{sec:CLR_MP}

I present for simplicity the spirit of the numerical methods developed in~\cite{PlStWa21} for general nondegenerate stochastic dynamics on~$\mathcal{D} = \mathbb{T}^d$, but the approach can be extended to degenerate dynamics such as Langevin dynamics, as discussed in~\cite{PlStWa22}. The reference dynamics
\[
dX^0_t = b(X^0_t) \, dt + \sigma(X^0_t) \, dW_t
\]
is perturbed as
\[
dX^\eta_t = \left( b(X^\eta_t) + \eta F(X^\eta_t) \right) dt + \sigma(X^\eta_t) \, dW_t.
\]
Assume that, for any~$\eta \in \mathbb{R}$, there exists a unique invariant probability measure~$\nu_\eta$ for this dynamics, which is the case for instance when~$b,\sigma$ are smooth, $\sigma \sigma^\top$ is everywhere positive definite, and the configuration space is compact. Of course, much milder conditions can be considered.

The key equality on which the numerical analysis is based is the following alternative expression for the transport coefficient, for an observable with average~0 with respect to the invariant probability measure~$\nu_0$ of the reference dynamics: 
\begin{equation}
  \label{eq:alpha_MP}
  \alpha = \lim_{\eta \to 0} \frac{\nu_\eta(R)}{\eta} = \lim_{t \to \infty}\mathbb{E}_0\left\{ \left(\frac{1}{t}\int_0^t R(X_s^0) \, ds\right) Z_t \right\},
\end{equation}
with
\[
Z_t = \int_0^t U(X_s^0)^\top dW_s, \qquad \sigma U = F.
\]
The above formula can be motivated by rewriting expectations for the perturbed dynamics as expectations for the reference dynamics weighted by a Girsanov factor as
\[
\begin{aligned}
  & \mathbb{E}_\eta\left[ \frac1t \int_0^t R(X_s^\eta)\, ds \right] \\
  & \quad =  \mathbb{E}_0\left[ \left(\frac1t \int_0^t R(X_s^0)\, ds\right) \exp\left(\eta\int_0^t U(X_s^0)^\top dW_s - \frac{\eta^2}{2}\int_0^t \left| U(X_s^0) \right|^2 ds  \right )\right],
  \end{aligned}
\]
then formally linearizing the resulting formula with respect to the small parameter~$\eta$ (in which case the exponential term in the above equality is replaced by~$1+\eta Z_t$ at dominant order in~$\eta$), and finally passing to the longtime limit. 

To rigorously prove the consistency of the approach, one introduces the generator~$\mathcal{L} + \eta \widetilde{\mathcal{L}}$ of the perturbed dynamics, and the Poisson equation $-\mathcal{L}\mathscr{R} = R$ (assumed to be well posed). As in~\eqref{eq:A_eta_t_Ito} and~\eqref{eq:GK_Ito}, the time integral of~$R$ can be rewritten as a martingale, up to remainder terms:
\[
\int_0^t R(X_s^0) \, ds = M_t + \mathscr{R}(X_0^0)-\mathscr{R}(X_t^0), \qquad M_t = \int_0^t \nabla \mathscr{R}(X_s)^\top \sigma(X_s^0) \,dW_s.
\]
It then remains to use It\^o's isometry to write $t^{-1} \mathbb{E}\left(M_t Z_t \right)$ as 
\[
\begin{aligned}
  & \frac1t \int_0^t \mathbb{E}\left[U(X_s^0)^\top \sigma(X_s^0)^\top \nabla \mathscr{R}(X_s^0) \right] ds \\
  & \qquad \xrightarrow[t\to+\infty]{} \int_\mathcal{D} F^\top \nabla \mathscr{R} \, d\nu_0 =  \int_\mathcal{D} \widetilde{\mathcal{L}} \mathscr{R} \, d\nu_0  = \int_\mathcal{D} \mathscr{R} \left(\widetilde{\mathcal{L}}^*\mathbf{1} \right) d\nu_0 = \alpha,
\end{aligned}
\]
where adjoints are considered on the Hilbert space~$L^2(\nu_0)$, and where the last equality follows by computations similar to the ones used to establish~\eqref{eq:LR_formula}. Similar manipulations allow to show that the variance of the estimator is uniformly bounded in time:
\[
\forall t > 0, \qquad \mathrm{Var}\left\{ \left(\frac{1}{t}\int_0^t R(X_s^0) \, ds\right) Z_t \right\} \leq C.
\]

From a numerical viewpoint, the limit~\eqref{eq:alpha_MP} suggests to consider the following discrete estimator (slightly idealized as the response function is centered according to the invariant probability measure of the numerical scheme): 
\[
\mathcal{M}_{{\Delta t}, N_{\rm iter}}^{[1]} = \frac{1}{N_{\rm iter}}\sum_{n=0}^{N_{\rm iter}-1} \left( R(X^n) - \mathbb{E}_{{\Delta t}}(R)\right ) Z^{N_{\rm iter}},
\]
where
\[
Z^{N_{\rm iter}} = \sum_{n=0}^{N_{\rm iter}-1} \left(\sigma(X^n)^{-1}F(X^n)\right)^\top G^n.
\]
It can then be shown that there exist~$C_1,C_2,C_3 \in \mathbb{R}_+$ such that, for~$\Delta t$ sufficiently small, 
\[
\begin{aligned}
  \left|\mathbb{E}_{{\Delta t}}\left\{\mathcal{M}_{{\Delta t}, N_{\rm iter}}^{[1]}\right\} - \alpha \right| & \leq C_1 \left({\Delta t} + \frac{1}{\sqrt{N_{\rm iter}{\Delta t}}}\right), \\
  \mathrm{Var}_{{\Delta t}}\left\{\mathcal{M}_{{\Delta t}, N_{\rm iter}}^{[1]}\right\} & \leq C_2 + C_3\left({\Delta t} + \frac{1}{N_{\rm iter}{\Delta t}}\right).
\end{aligned}
\]
These estimates show that the bias of the estimator is of order~$\Delta t$ (although it can be reduced to~$\Delta t^2$ upon considering schemes of weak order~2 and correcting the discrete martingale term, see~\cite{PlStWa22}), and that the bias arising from the finite time integration scales as the inverse of the square root of the physical time, instead of scaling as the inverse of the physical time as estimates such as~\eqref{eq:FT_bias}. This is due to the extra factor~$Z^{N_{\rm iter}}$ multiplying the time average. 

\subsection{Current perspectives on better estimating transport coefficients}
\label{sec:current_perspectives}

Let me conclude this review by listing various alternative strategies to the usual approaches discussed in Section~\ref{sec:error}: 
\begin{itemize}
\item Some variance reduction can be obtained in nonequilibrium molecular dynamics by considering control variate approaches. There are various realizations of this idea. As discussed in~\cite{RoSt19} (see also~\cite{MaRo20}), is possible for instance to consider a response function~$R+(\mathcal{L}+\eta\widetilde{\mathcal{L}})\Phi$ instead of~$R$. The idea is that the modified response function still has the same average as~$R$ under the invariant probability measure of the nonequilibrium dynamics, but can have a much smaller asymptotic variance when~$\Phi$ is well chosen. Dynamical versions of the control variate idea can also be considered~\cite{CiJa75,PaStVa22}. 
\item It may be beneficial to use a coupling between the perturbed dynamics~$X_t^\eta$ and the reference one~$X_t^0$, for instance sticky coupling~\cite{EbZi19,DuEbEnGuMo21}. This idea is currently explored in~\cite{BoDaEbSt22}.
\item One can also rely on tangent dynamics~\cite{AsJoLeRo18}, which consist in formally passing to the limit~$\eta \to 0$ in two dynamics which are synchronously coupled. The resulting evolution equations are the reference dynamics on~$X_t^0$, and a random ordinary evolution for the tangent vector~$T_t = \lim_{\eta \to 0} (X_t^\eta-X_t^0)/\eta$. A Green--Kubo type formula can be written in terms of~$X_t^0$ and~$T_t$.
\item While physically intuitive external forcings are most often considered in nonequilibrium molecular dynamics simulations, it is in fact possible to optimize the nature of the forcing. The key remark here is that infinitely many forcings give the same linear response as the physical one. In essence, forcings can differ by any operator which preserves the invariant measure, as this leaves the conjugate response function~\eqref{eq:conjugate_response} unchanged. Such forcings are called ``synthetic forcings'' in~\cite{EvMo08}. From a mathematical perspective, this suggests to optimize the extra forcing added to the system, in order for instance to increase the linear response regime; see~\cite{SpSt22}. 
\item Large deviation techniques have also been used to estimate cumulant generating functions of time averages of a response function~$R$. In particular, second order cumulants, which correspond to asymptotic variances or Green--Kubo formulas, are obtained by a quadratic approximation around~0. The idea here is that, in order to have a good polynomial fit around~0, it may be beneficial to determine the cumulant generating function away from~0; see~\cite{LiGaPo21}.
\end{itemize}
Of course, these approaches can be combined. One can imagine for instance considering the tangent dynamics associated with a coupling different than the synchronous coupling; or using control variate methods in conjunction with synthetic forcings. There are also other approaches which are yet too prospective to be mentioned at this stage...

The success or failure of possible alternative methods should be asserted in any case based on two indicators:
\begin{enumerate}[(i)]
\item theoretical estimates which quantify the numerical errors, in terms of variance and bias, as a function of the parameters of the numerical method (value of the forcing magnitude~$\eta$, timestep~$\Delta t$, number of iterations~$N_{\rm iter}$, etc);
\item numerical simulations on representative systems of interest, such as Lennard--Jones fluids or atom chains (the thermal conductivity of the latter systems being rather challenging to properly estimate as the system sizes are increased, see for instance~\cite{LeLiPo03,Dh08,IuLeLiPoPo19}).
\end{enumerate}
I modestly hope that this review will trigger interest in the community to tackle these issues!

\begin{acknowledgement}
  The author thanks Arnaud Guyader for useful and relevant feedback on the slides presented during the conference, as well as No\'e Blassel, Régis Santet and Renato Spacek for carefully re-reading a preliminary version of this review. The remarks and comments of the anonymous referee also allowed to improve the current version of the manuscript. This project has received funding from the European Research Council (ERC) under the European Union's Horizon 2020 research and innovation programme (project EMC2, grant agreement No 810367), and from the Agence Nationale de la Recherche, under grants ANR-19-CE40-0010-01 (QuAMProcs) and ANR-21-CE40-0006 (SINEQ). 
\end{acknowledgement}

%-- biblio -----
%please use a prefix corresponding to the surnames of the authors of your paper: E.g., if the paper is single-authored by Minnie Mouse, all labels should have the prefix "Mo". If the authors of a paper are Donald Duck and Bart Simpson, all labels should have the prefix "DuSi", and accordingly for more authors. This will help us a lot in avoiding clashes between the labels used in different papers.

\end{document}